\documentstyle[12pt, leqno,amssymb,amscd,graphics, epstopdf ,graphicx,geometry]{amsart}

\newtheorem{thm}{Theorem}[section]
\newtheorem{lemma}[thm]{Lemma}

\theoremstyle{definition}
\newtheorem{dfn}[thm]{Definition}
\theoremstyle{remark}

\begin{document}

\newcommand{\ct}{\cite}
\newcommand{\pr}{\protect\ref}
\newcommand{\su}{\subseteq}
\newcommand{\pa}{{\partial}}
\newcommand{\e}{\epsilon}
\newcommand{\es}{{\varnothing}}
\newcommand{\D}{{\mathcal D}}

\newcommand{\lk}{{\mathcal{L}}}
\newcommand{\G}{{\mathbb G}}
\newcommand{\Z}{{\mathbb Z}}
\newcommand{\R}{{\mathbb R}}
\newcommand{\SP}{{\mathbb S}}
\newcommand{\I}{{\mathrm{Id}}}
\newcommand{\dg}{{\mathcal{D}}}

\newcounter{numb}

\title{Complexity of relations in the braid group}

\author{Joel Hass}
\address{Department of Mathematics, University of California, Davis, CA 95616, U.S.A.} 
\email{hass@@math.ucdavis.edu}
\urladdr{}

\author{Arkadius Kalka}
\address{Department of Mathematics, Bar-Ilan University, 
Ramat-Gan 52900, Israel}
\email{Arkadius.Kalka@@ruhr-uni-bochum.de}
\urladdr{}

\author{Tahl Nowik}
\address{Department of Mathematics, Bar-Ilan University, 
Ramat-Gan 52900, Israel}
\email{tahl@@math.biu.ac.il}
\urladdr{}

\date{May 28, 2009}

\begin{abstract}
We show that for any given $n$, there exists a sequence of words $\{ a_k \}_{k \geq 1}$ in the generators
$\sigma_1,\dots,\sigma_{n-1}$ of the braid group $B_n$, representing the identity element of
$B_n$, such that the number of braid relations of the form
$\sigma_i\sigma_{i+1}\sigma_i = \sigma_{i+1}\sigma_i\sigma_{i+1}$ needed to pass from $a_k$ to the empty
word is quadratic with respect to the length of $a_k$.
\end{abstract}

\maketitle

\section{Introduction}\label{int}
We introduce in this paper a new method to understand the complexity of  word problems in certain groups.
This problem, and the closely associated problem of finding a Dehn function for a group, have been a focus of study in combinatorial group theory. 
We then apply our method to obtain new information  about the complexity of certain words in the braid group.

Suppose that  $w$ is a word in a finitely presented group $G$ which represents the identity element and  
is a product of $k$ generators. Then there is a sequence of  $f(k)$ relators that can be applied to $w$ to 
produce the empty word. A Dehn function gives an upper bound for $f(k)$ in terms of $k$, and there has been extensive
study of these functions \cite{BRS},\cite{BORS},\cite{bb},\cite{ECHLPT},\cite{Ge},\cite{G}.

While the Dehn function indicates the worst case for words of length $k$, it doesn't give information about 
particular words. We present in this paper a technique which, in certain cases, 
allows the counting of the number of relations that need to be applied to transform a particular word $w$ 
to the empty word.
An invariant is assigned to a word and the effect of a relator on the invariant is determined.
The invariant can then give lower bounds on the number of relators needed to connect
two words.

We apply this to study certain words in the braid group.  The braid group on four or more strands has a noncyclic free abelian subgroup, so its Dehn function is at least quadratic.   The abelian subgroups arise from crossings between pairs of non-adjacent strands, which lead to commutativity relations. 
We show that the braid groups have non-linear Dehn function for deeper reasons.  The number of relations required to trivialize certain words  is quadratic even without counting the commutativity relations.

To get these results we consider certain invariants assigned to words in the braid group.  
These invariants are obtained from invariants of knot diagrams that were introduced in \cite{hn}.  
We compute how these invariants change when a relator in the braid group is applied to a word, 
showing that they are unchanged by the commutativity relations and by the canceling of a generator with its inverse,
and change by $\pm 1$ under the other relations.  
We then produce a sequence of words $\{ a_k \} $
which represent the identity element and for which the invariant is 
quadratic with respect to the length of $a_k$.  
This implies that at least a quadratic number of relations are needed to trivialize $a_k$, 
even with commutativity relations not counted.

We note that Thurston showed that braid groups are automatic, and that the the word problem in $B_n$ is solvable in $O(k^2 n \log n)$  steps \cite{ECHLPT}.  Therefore there exists also a quadratic upper bound for these words.

\section{Statement of result}\label{st}

It is known that the braid group $B_n$ on $n$ strands has the following presentation:
$$\langle \sigma_1,\dots,\sigma_{n-1} \ | 
\ \sigma_i\sigma_{i+1}\sigma_i = \sigma_{i+1}\sigma_i\sigma_{i+1} \ \hbox{for} \ i=1,\dots,n-2, \ \hbox{and} \ 
\sigma_i\sigma_j = \sigma_j\sigma_i \ \hbox{for} \ |i-j| \geq 2 \rangle$$
where we take $\sigma_i$ to be a positive crossing between the $i$th and $(i+1)$th strand, as appears 
in Figure \pr{f1}.

\begin{figure}[h]
\scalebox{0.8}{\includegraphics{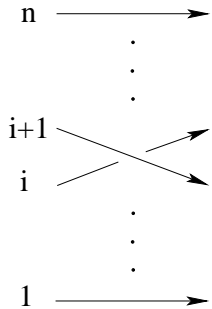}}
\caption{The generator $\sigma_i$.}\label{f1}
\end{figure}

If $a$ and $b$ are two words in the generators $\sigma_1,\dots,\sigma_{n-1}$ representing the same element in $B_n$, 
then there is a finite sequence of modifications that brings us from $a$ to $b$, where each modification
is of one of the following three types: 
\begin{enumerate}
\item A subword of the form $\sigma_i\sigma_{i+1}\sigma_i$ is replaced by $\sigma_{i+1}\sigma_i\sigma_{i+1}$ 
or vice versa. 
\item A subword of the form $\sigma_i\sigma_j$ is replaced by $\sigma_j\sigma_i$ for $|i-j| \geq 2$.
\item A subword of the form $\sigma_i\sigma_i^{-1}$ or $\sigma_i^{-1}\sigma_i$ is replaced by the empty word, 
or vice versa.
\end{enumerate}

We will be interested in counting the number of modifications of type (1) required to pass from $a$ to $b$.
So, we define:

\begin{dfn}\label{d}
Given words $a$ and $b$ in the generators $\sigma_1,\dots,\sigma_{n-1}$ of $B_n$, 
representing the same element in $B_n$, we define the distance $d(a,b)$ between $a$ and $b$,
to be the minimal number of modifications of type (1) that may appear in a sequence of modifications
passing from $a$ to $b$.
\end{dfn}

We note that we could have included as basic modifications also those interchanging between  
$\sigma_i^{-1}\sigma_{i+1}^{-1}\sigma_i^{-1}$ and $\sigma_{i+1}^{-1}\sigma_i^{-1}\sigma_{i+1}^{-1}$
and between $\sigma_i^{-1}\sigma_j^{-1}$ and $\sigma_j^{-1}\sigma_i^{-1}$ for $|i-j| \geq 2$.
But this would not have changed the definition of $d(a,b)$, since each such modification can be achieved 
by a sequence of modifications which includes precisely one modification of the same type with positive exponents, 
and all other modifications of type (3).

For $n<m$ there is a natural embedding $B_n \hookrightarrow B_m$ given by adding $m-n$ straight strands to the 
given braid. In terms of the generators this is given simply by 
$\sigma_i \mapsto \sigma_i$, $i=1,\dots, n-1$. We now show that the distance $d(a,b)$
is preserved under this embedding.

\begin{lemma}\label{emb}
Let $a,b$ be two words in $\sigma_1,\dots,\sigma_{n-1}$ representing the same element in $B_n$. 
Then the distance $d(a,b)$ is the same whether $a,b$ are interpreted in $B_n$, or in $B_m$ for $m>n$.
\end{lemma}

\begin{pf}
Denote by $d_n,d_m$ the distances in $B_n,B_m$ respectively.
Any sequence of modifications in $B_n$ from $a$ to $b$ defines the same sequence in $B_m$, 
so $d_n(a,b) \geq d_m(a,b)$. 
On the other hand, given a sequence in $B_m$ from $a$ to $b$, 
look at the braid diagram of each word in $B_m$ along the way, and delete from it
its $m-n$ last strands. This produces a sequence in $B_n$ from $a$ to $b$. 
Some of the modifications in $B_m$ induce the corresponding modifications in $B_n$, whereas 
others induce no modification at all. We get $d_m(a,b) \geq d_n(a,b)$, so together $d_m(a,b) = d_n(a,b)$. 
\end{pf}

In this work we prove the following:

\begin{thm}\label{main}
For any $n$ there exists a sequence of words $\{ a_k \}_{k \geq 1}$ in the generators
$\sigma_1,\dots,\sigma_{n-1}$ of $B_n$, representing the identity element of
$B_n$, such that the length of the word $a_k$ 
is $12k$, and $d(a_k,E) = 4k^2$ where $E$ is the empty word. 
(That is, the minimal number of modifications of the form
$\sigma_i\sigma_{i+1}\sigma_i = \sigma_{i+1}\sigma_i\sigma_{i+1}$ required to pass from $a_k$ to $E$
is $4k^2$.) 
\end{thm}

\section{Proof of Theorem \pr{main}}\label{prf}

By Lemma \pr{emb} it is enough to prove Theorem \pr{main} for $B_3$.
We will prove it as follows. For any word $w$ in $B_3$ representing the identity element, we define
the knot diagram $\widehat{w}$ to be the closure of the braid diagram $w \sigma_2\sigma_1$.
(The purpose of adding $\sigma_2\sigma_1$ is that
a knot will be obtained, rather than a three component link.)
Our braid modifications on $w$ induce Reidemeister moves on the knot diagram $\widehat{w}$. We distinguish 
four types of Reidemeister moves, RI, matched RII, unmatched RII, and RIII,
as appears in Figure \pr{f2}.
We will use the invariant of knot diagrams introduced in \ct{hn} 
to analyze our braid modifications via the Reidemeister moves that they induce.

\begin{figure}[h]
\scalebox{0.8}{\includegraphics{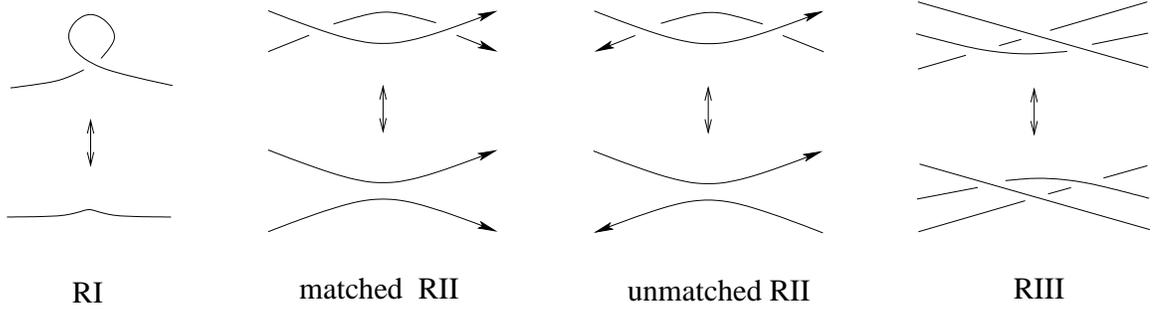}}
\caption{Reidemeister moves.}\label{f2}
\end{figure}

We first recall the definition of the invariant of \ct{hn}.
We denote the set of all knot diagrams by $\dg$ and
the set of all two component links in $\R^3$ by $\lk$. 
Given a knot diagram $D \in \dg$ and a crossing $a$ in $D$, define $D^a \in \lk$
to be the two component link obtained by smoothing the crossing $a$.
Given a knot diagram $D$, let   $D_+$ denote the set of all positive
crossings in $D$ and  $D_-$ the set of all negative crossings.
Given an invariant of two component links $\phi:\lk \to S$ where $S$ is any set,
let $\G_S$ be the free abelian group with basis $\{X_s, Y_s\}_{s \in S}$. We then define 
the invariant $I_\phi : \dg \to \G_S$ to be
$$I_\phi(D) = \sum_{a \in D_+} X_{\phi(D^a)}  + \sum_{a \in D_-} Y_{\phi(D^a)}.$$

In the present work we take  $\phi$ to be the linking number, $lk:\lk\to\Z$,  giving the invariant
$I_{lk} : \dg \to \G_\Z$. 
In \ct{hn} it is shown that the change in the value of $I_{lk}$ 
resulting from the various types of Reidemeister moves is as follows:
\begin{itemize}
\item For an RI move: $X_0$ or $Y_0$. 
\item For a matched RII move: $X_j + Y_{j+1}$.
\item For an unmatched RII move: $X_j + Y_j$.
\item For an RIII move: $X_j - X_{j+1}$ or $Y_j - Y_{j+1}$.  
\end{itemize}

Our modification of braid words of type (1) produces a Reidemeister move of type RIII,
our modification of type (2) produces no change in the knot diagram, and our modification of 
type (3) produces a matched RII move. 
Let $g:\G_\Z \to \Z$ be the homomorphism given by
$g(X_j)=|j|$ and $g(Y_j)=-|j-1|$. Then  
$g(X_j - X_{j+1})=\pm 1$,   $g(Y_j - Y_{j+1})= \pm 1$, and $g(X_j + Y_{j+1})=0$. 

We now define an invariant $G$ of braid words representing the identity element in $B_3$, as follows:  
$G(w) = g\circ I_{lk} (\widehat{w})$.
It follows from the observations of the previous paragraph that $G(w)$ 
changes by $\pm 1$ whenever a modification of type 1 takes place, and remains unchanged whenever a 
modification of type 2 or 3 takes place. It follows that for any 
two words $a,b$ representing the identity element in $B_3$,
$d(a,b) \geq |G(a) - G(b)|$.

\begin{figure}[t]
\scalebox{0.8}{\includegraphics{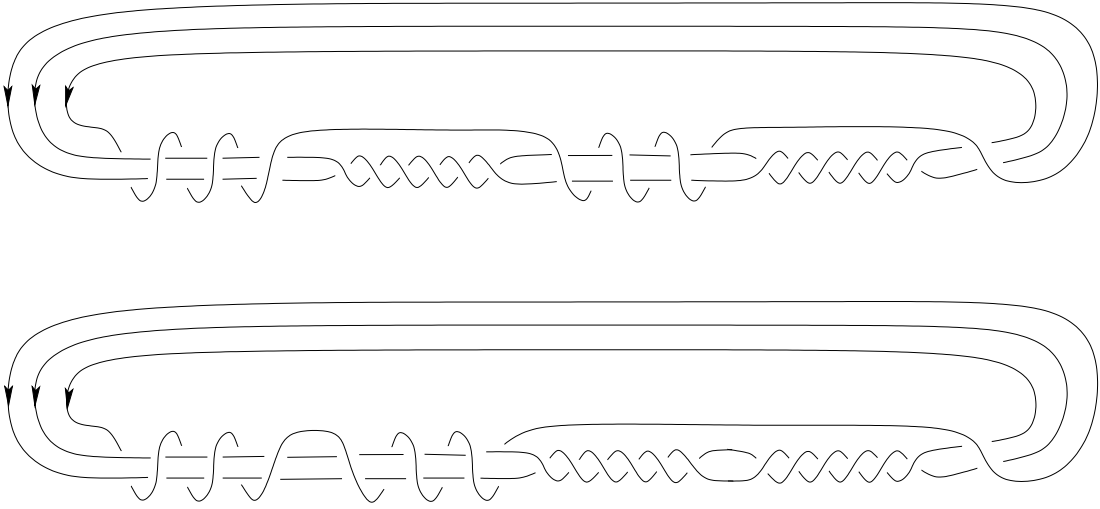}}
\caption{The knot diagrams $\widehat{a_k}$ and $\widehat{b_k}$ for $k=3$.}\label{f3}
\end{figure}

We now construct our sequence of words $a_k$, together with a second sequence of words $b_k$. 
We will show that $d(a_k,b_k)=4k^2$, and it will be evident that $d(b_k,E)=0$. Since $d$ is a pseudometric
it follows that $d(a_k,E)=d(a_k,b_k)=4k^2$.
We construct our two words $a_k,b_k$ in the generators $\sigma_1,\sigma_2$ of $B_3$ as follows.
Our building blocks will be the two words $V=\sigma_1^2$ and  $W=\sigma_2\sigma_1^2\sigma_2$, and we define
$a_k = W^{-k} V^k W^k V^{-k}$ and $b_k = W^{-k} W^k V^k V^{-k}$. Note that the length of $a_k$ and $b_k$ 
in terms of the generators $\sigma_1,\sigma_2$ is $12k$.
The knot diagrams $\widehat{a_k}$ and $\widehat{b_k}$ are depicted in Figure \pr{f3}.
(Recall that by definition they are the closures of the braid diagrams  $a_k\sigma_2\sigma_1$ and 
$b_k\sigma_2\sigma_1$ respectively.) 
From the figure, it is clear that the portion $V^k$ in $a_k$ can slide through the portion $W^k$
to obtain $b_k$. So indeed $a_k$ and $b_k$ both represent the identity element in $B_3$. 
This sliding process requires precisely $4k^2$ applications of the relation
$\sigma_1\sigma_2\sigma_1=\sigma_2\sigma_1\sigma_2$ and so $d(a_k,b_k) \leq 4k^2$.

\begin{figure}[t]
\scalebox{0.8}{\includegraphics{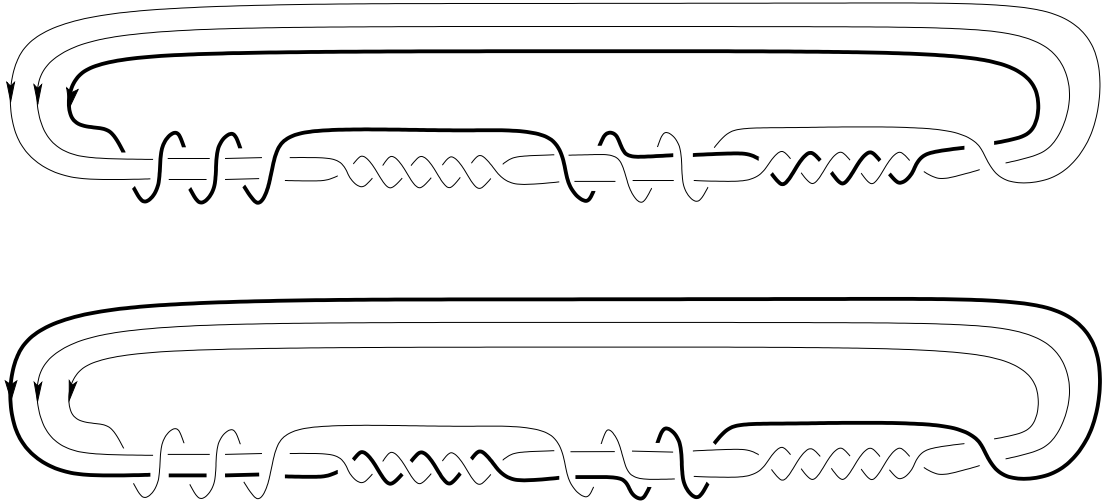}}
\caption{Two smoothings of $\widehat{a_k}$.}\label{f4}
\end{figure}

As explained above, a lower bound for $d(a_k,b_k)$ may be obtained by computing $|G(a_k)-G(b_k)|$.
For this purpose we first compute $I_{lk}(\widehat{a_k}) - I_{lk}(\widehat{b_k})$ 
and then apply $g$ to this element of $\G_\Z$.
We note that whenever we smooth a given crossing in the initial portion
$W^{-k}$ of $\widehat{a_k}$, or in its concluding portion $V^{-k}\sigma_2\sigma_1$, then
the link obtained is isotopic to that obtained from $\widehat{b_k}$ by the corresponding smoothing,
this isotopy being given by sliding the $V^k$ portion through the $W^k$ portion.
Furthermore, any smoothing in the $V^k$ portion also gives isotopic links for $\widehat{a_k}$ and $\widehat{b_k}$
since after such a smoothing this portion may still be slid through the $W^k$ portion.
Each such crossing whose smoothing produces isotopic links for $\widehat{a_k}$ and $\widehat{b_k}$,
contributes the same term to $I_{lk}(\widehat{a_k})$ and $I_{lk}(\widehat{b_k})$, which cancels in the difference
$I_{lk}(\widehat{a_k}) - I_{lk}(\widehat{b_k})$.
So to compute $I_{lk}(\widehat{a_k}) - I_{lk}(\widehat{b_k})$ we must only consider crossings in the $W^k$ portion.
Two such representative smoothings for $\widehat{a_k}$ are depicted in Figure \pr{f4},
and the corresponding two smoothings for $\widehat{b_k}$ are depicted in Figure \pr{f5}.

Each $W$ in the $W^k$ portion has four crossings. The top of Figure \pr{f4} demonstrates the smoothing of the first 
crossing in, say, the $i$th $W$. One may see from the figure 
that the linking number of the link obtained is 
$\frac{1}{2} \bigl( -4k+4(i-1)+1+2(k-i)-2k+1 \bigr) =i-2k-1$. Since the crossing itself is positive, this
contributes $X_{i-2k-1}$ to $I_{lk}(\widehat{a_k})$. 
The bottom of Figure \pr{f4} demonstrates the smoothing of the second crossing of the same $W$.
One may see that the linking number here is
$\frac{1}{2} \bigl( -2k +2k +2(i-1) + 2 + 4(k-i) + 2 \bigr) = 2k+1-i$, and so this crossing, being itself positive, 
contributes $X_{2k+1-i}$.
The third crossing produces a knot diagram identical to that of the second crossing, so also contributes $X_{2k+1-i}$.
The fourth crossing produces a knot diagram resembling that of the first crossing, but with two additional 
positive crossing between the two components of the link, and  so 
its linking number is $i-2k$. (In fact, for $i<k$, the smoothing of the fourth crossing of the $i$th $W$
produces an identical knot diagram to that of the first crossing of the $(i+1)$th $W$.)
So, the fourth crossing contributes $X_{i-2k}$.

\begin{figure}[t]
\scalebox{0.8}{\includegraphics{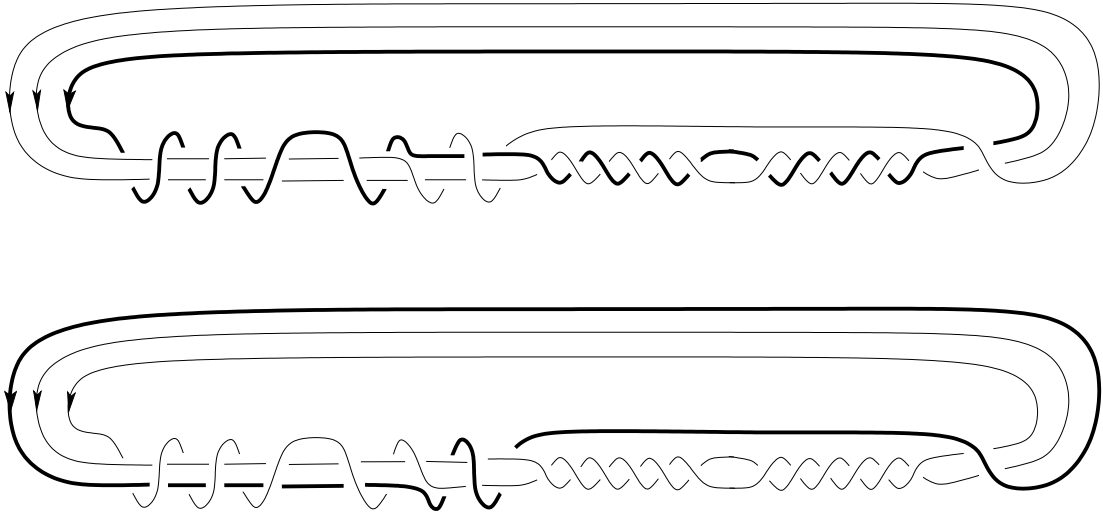}}
\caption{Two smoothings of $\widehat{b_k}$.}\label{f5}
\end{figure}

Similar smoothings for $\widehat{b_k}$ appear in Figure \pr{f5}.
One may verify that the contribution of the four crossings of the $i$th $W$ are in this case
$X_{i-k-1}$, $X_{k-i+1}$ ,$X_{k-i+1}$, and $X_{i-k}$, respectively. 

Summing over $i$ gives:
$$ I_{lk}(\widehat{a_k}) - I_{lk}(\widehat{b_k})
=\sum _{i=1}^{k} (X_{i-2k-1}+2X_{2k+1-i}+X_{i-2k}-X_{i-k-1}-2X_{k-i+1}-X_{i-k}). $$
Applying $g$ we obtain $G(a_k)-G(b_k) = g \bigl( I_{lk}(\widehat{a_k}) - I_{lk}(\widehat{b_k}) \bigr) =4k^2$, 
which proves our lower bound, and so $d(a_k,b_k)=4k^2$. 
As explained above, since clearly $d(b_k,E)=0$, it follows that also $d(a_k,E)=4k^2$. 
This completes the proof of Theorem \pr{main}.

We remark that $d(a_k,b_k)=4k^2$ implies that also $d(V^kW^k,W^kV^k)=4k^2$, though these words
do not represent the identity.

\end{document}